# Approximations of Time-Dependent Nonlinear Partial Differential Equations using Galerkin Optimal Auxiliary Function Method

Nilormy Gupta Trisha and Md. Shafiqul Islam*

Department of Applied Mathematics, University of Dhaka, Dhaka 1000, Bangladesh.

Email : nilormy1997@gmail.com, *Corresponding author and Email : mdshafiqul@du.ac.bd

## Abstract

The purpose of this research work is to employ the Optimal Auxiliary Function Method (OAFM) for obtaining numerical approximations of time-dependent nonlinear partial differential equations (PDEs) that arise in many disciplines of science and engineering. The initial and first approximations of parabolic nonlinear PDEs associated with initial conditions have been generated by utilizing this method. Then the Galerkin method is applied to estimate the coefficients that remain unknown. Finally, the values of the coefficients generated by the Galerkin method have been inserted into the first approximation. In each example, all numerical computations and corresponding absolute errors are provided in schematic and tabular representations. The rate of convergence attained by the proposed method is depicted in tabular form.

## 1 Introduction

Nonlinear parabolic partial differential equations (PDEs) are used in simulating a wide variety of physical phenomena in the technological and scientific realms, from turbulence to particle dispersion, as well as in establishing the valuation of a wide range of derivative financial instruments. They are utilized in the endeavor of providing an explanation for a wide range of occurrences, including liquid filtration, sound, heat, diffusion, chemical reactions, fluid dynamics, environmental contamination, and many more.
Nonlinear parabolic PDEs have been numerically studied using the well-known Adomian Decomposition Method [1]. The adaptive grid Haar wavelet collocation method [2] has been used to get quantitative solutions to these types of equations. Lang [3] has authored a book on the applicability of multilevel solutions to the parabolic PDE system. For fully nonlinear PDEs, Arash Fahim et al. [4] developed a method that combines Monte Carlo with a finite difference scheme. JW Wang et al. [5] have published an article utilizing a fuzzy control approach to address the complexities of nonlinear parabolic PDE systems. Several authors have put together a comprehensive reference work on these nonlinear PDEs. This book [6] covers a wide range of approaches and various aspects of these equations. Various eminent authors have published numerous books on these nonlinear PDEs [7, 8, 9, 10]. They have thoroughly examined these types of equations and analyzed their relevance to real-world issues. In order to address some particularly challenging instances of parabolic PDEs, the field of wavelet analysis has recently emerged as a significant mathematical tool [11]. Ekren has provided the viscosity solutions for fully nonlinear parabolic spatial PDEs [12]. In order to deal with the nonlinearity of parabolic PDEs, Mironchenko et al. [13] have utilized the monotony-control system. To estimate solutions to a set of nonlinear parabolic PDEs, Izadi et al. [14] have proposed an innovative hybrid spectral collocation methodology. In the research study cited in [15], an expansion of the Taylor series has been presented as a method for numerically solving parabolic PDEs. Many renowned authors [16, 17, 18, 19] have approximated parabolic PDEs associated with different boundary conditions with the help of the Galerkin Finite Element Method. Kamrujjaman et al. [20] have applied the finite difference method to renowned nonlinear PDEs. Alam et al. [21] have provided the approximate solutions of different parabolic PDEs—heat and wave equations.

Each of the above approaches offers benefits and drawbacks. Likewise, we generalize a recently established methodology called the optimal auxiliary function method (OAFM) to the context of PDEs. The pioneers of this approach were Marinca & Marinca [22]. For obtaining an analytical solution of the fluid's thin layer in the cylinder's vertical orientation, they used this methodology. This method offers a means of controlling the convergence of numerical solutions with the assistance of convergence-control parameters in order to achieve the desired level of accuracy. Asserting such a command allows for controlled execution. It is a straightforward, convergent, moderate, and explicit method for obtaining nonlinear approximations. A further application of this approach was made by Marinca et al. [23] to resolve the nonlinear Blasius issue. In addition to this, it has been used to analyze the nonlinear vibrations of a pendulum that has been folded around two cylinders [24]. This method has been offered as an approximate analytical solution for the nonlinear boundary issue of viscous flow that is created by a stretched surface with partial slippage [25]. The thin layer of a third-grade fluid has been modeled using the optimal auxiliary function method so that it can be simulated on a moving belt [26]. It has also been demonstrated that the technique is superior to other approaches in terms of its effectiveness in addressing challenges brought on by misalignment [27]. Recently, Laiq Zada has utilized this method in order to create approximate-analytical solutions to partial differential equations and generalized modified b-equations, as referenced in [28] and [29]. In reference [30], the OAFM is generalized to the realm of partial differential equations and is employed to approximately solve KdV equations. Ullah et al. [31] have just recently utilized this strategy in order to evaluate fractional KdV equations. In reference [32], the proposed method has been employed for the steady nanofluids. The method has recently been applied in biological modeling [33]. In the research studies cited in [34] and [35], a thorough study into the electronic and physical modifications of a low-power PMS generator has been carried out with the assistance of the OAFM.

The application of the Galerkin Weighted Residual Method (GWRM) dates back centuries, even before the advent of computers. The method is acknowledged to be one of the best and most widely applied methods. In the book cited in [36], Lewis and Ward have provided a detailed description of the procedure. Hossan et al. [37] have successfully implemented this strategy on the well-known Black-Scholes model. In their analysis of the Fredholm equations, Shirin et al. [38] have used the Galerkin method in conjunction with other specialized polynomials. The method has been applied to the boundary value problems in the research cited in [39]. It has also been employed to numerically calculate the eigenvalues of the Sturm-Liouville problem [40]. The technique has had extensive application in issues involving metal beams and polygonal ducts with rounded corners. [41, 42].

Inspired by all previous research, our proposed methodology has been deployed to solve some renowned nonlinear parabolic PDEs that are associated with corresponding initial conditions. So far, the method has not required any complicated computations due to the variety of boundary conditions considered.

As far as we know, this approach to solving these kinds of PDEs is not available. All approximate results are analyzed in light of the exact solutions of the parabolic PDEs that have been provided.
This article is split into four distinct sections. Section 2 provides a formal derivation of our suggested approach. In the third section, the approach's implications are shown while analyzing four problems characterized by high nonlinearity. Graphs of errors and numerical findings are included here as well. The fourth section contains some concluding remarks and a general discussion.

## 2 Mathematical Formulation

Let us start with a general PDE of the form [28]:

$$\Lambda[M(x,t)] + \Upsilon[M(x,t)] + g(x,t) = 0, \qquad x \in \mathcal{D} \tag{1}$$

subject to the boundary/initial conditions $\Omega\left[M, \dfrac{\partial M}{\partial t}\right] = 0$.
Here $\Lambda$ is considered the linear operator, and $\Upsilon$ is regarded as the nonlinear operator. Again we consider $g(x,t)$ and $M(x,t)$ as known function and unknown functions respectively. The domain of interest is $\mathcal{D}$.
Let the form of the approximate solution of Equation (1) be

$$\widetilde{M}(x,t) = M_0(x,t) + M_1(x,t,C_i), i = 1,2,3...,p \tag{2}$$

where $C_i$'s are currently unknowable $p$ parameters. Here $M_0(x,t)$ is the initial approximation and $M_1(x,t)$ is the first approximation. In this case, $p$ is a positive number chosen at random.
Substituting (2) in Equation (1) we get,

$$\begin{aligned}
&\Lambda[\widetilde{M}(x,t)] + \Upsilon[\widetilde{M}(x,t)] + g(x,t) = 0\\
\text{or,}&\Lambda[M_0(x,t) + M_1(x,t,C_i)] + \Upsilon[M_0(x,t) + M_1(x,t,C_i)] + g(x,t) = 0\\
\text{or,}&\Lambda[M_0(x,t)] + \Lambda[M_1(x,t,C_i)] + \Upsilon[M_0(x,t)] + \Upsilon[M_1(x,t,C_i)] + g(x,t) = 0
\end{aligned} \tag{3}$$

Since the nonlinear operator is $\Lambda$, then we use the following linear equation to approximate the value of initial approximation $M_0(x,t)$ at the outset:

$$\Lambda[M_0(x,t)] + g(x,t) = 0, \qquad \Omega\left[M_0, \frac{\partial M_0}{\partial t}\right] = 0 \tag{4}$$

To estimate the first approximation of solution (2), we consider the second differential equation in the following form which also comprises the nonlinear operator $\Upsilon$,

$$\Lambda[M_1(x,t,C_i)] + \Upsilon[M_0(x,t)] + \Upsilon[M_1(x,t,C_i)] = 0, \qquad \Omega\left[M_1, \frac{\partial M_1}{\partial t}\right] = 0 \tag{5}$$

The nonlinear term in Equation (5) is expanded in the form as

$$\Upsilon[M_0(x,t)] + \Upsilon[M_1(x,t)] = \Upsilon[M_0(x,t)] + \sum_{\kappa\geq 1} \frac{M_1^{\kappa}(x,t,C_i)}{\kappa!} \Upsilon^{(\kappa)}(M_0(x,t)) \tag{6}$$

where $\kappa!$=1.2.3...$\kappa$ and $\Upsilon^{(\kappa)}$ stand for the differentiation of order $\kappa$ of the nonlinear operator $\Upsilon$. Instead of solving Equation (5), we want to find a way to circumvent the difficulties of resolving the nonlinear differential equation in (5) and hasten the first approximation's swift convergence and the implication of the solution $\widetilde{M}(x,t)$. Therefore, Equation (5) is shown to have an alternate expression

$$\Lambda[M_1(x,t,C_i)] + B_1(M_0(x,t),C_i)\Upsilon(M_0(x,t)) + B_2(M_0(x,t),C_j) = 0 \tag{7}$$

with

$$\Omega\left[M_1(x,t,C_i), \frac{\partial M_1(x,t,C_i)}{\partial t}\right] = 0 \tag{8}$$

$B_1$ and $B_2$ are presumed to be the random auxiliary functions that rely on the initial approximation $M_0(x,t)$ and unknown parameters $C_i$ and $C_j$, where $i = 1,2,...,s$ and $j = s+1, s+2, ..., p$ respectively. The positive integer $p$ and the auxiliary functions can be chosen in a wide variety of ways, providing us with a lot of flexibility. The auxiliary functions $B_1(M_0(x,t),C_i)$ and $B_2(M_0(x,t),C_j)$ are not unique; they rely on the initial approximation $M_0(x,t)$ or on the combinations of $M_0(x,t)$ and $\Upsilon[M_0(x,t)]$.
The unknown parameters' values $C_i$ and $C_j$ have been obtained by applying the Collocation method in the research article cited in [28]. But in this research article, we have employed the Galerkin method which is as follows:
First of all, the approximate solution $\widetilde{M}(x,t)$ can be redefined as :

$$\widetilde{M}(x,t) = M_0(x,t) + M_1(x,t) \tag{9}$$

$$= M_0(x,t) + \sum_{j=1}^{n} tC_j\phi_j(x) \tag{10}$$

where the functions $\phi_i(x)$ are called the coordinate functions. The **residual function** of Equation (1) can be written as,

$$R(x,t) = \Lambda[\widetilde{M}(x,t)] + \Upsilon[\widetilde{M}(x,t)] + g(x,t) \tag{11}$$

Then the residual function $R(x,t)$ by the coordinate functions and set the residual equation as

$$\int_{\mathcal{D}} R(x,t)\phi_i(x)dx = 0$$
$$\text{or} \int_{\mathcal{D}} \Big[\Lambda[\widetilde{M}(x,t)] + \Upsilon[\widetilde{M}(x,t)] + g(x,t)\Big]\phi_i(x)dx = 0 \tag{12}$$

Then we substitute the redefined solution (12) in the residual equation (12). It results in the following equation,

$$\int_{\mathcal{D}} \left[\Lambda\Big[M_0(x,t) + \sum_{j=1}^{n} tC_j\phi_j(x)\Big] + \Upsilon\Big[M_0(x,t) + \sum_{j=1}^{n} tC_j\phi_j(x)\Big] + g(x,t)\right]\phi_i(x)dx = 0$$
$$\text{or} \quad \sum_{j=1}^{n} C_jK_{ij} = F_i \tag{13}$$

where

$$K_{ij} = \int_{\mathcal{D}} \Big[\Lambda\big[t\phi_j(x)\big] + \Upsilon\big[t\phi_j(x)\big]\Big]\phi_i(x)dx \tag{14}$$

$$F_i = \int_{\mathcal{D}} \Big[\Lambda[M_0(x,t)] + \Upsilon[M_0(x,t)] + g(x,t)\Big]\phi_i(x)dx \tag{15}$$

The linear part of Equation (13) has been solved for the initial solution. Then the initial solutions are used for determining the coefficients $C_i$ and $C_j$. The method has been elaborated for the particular test problems in the following section.
With these known parameters, the approximated solution $\widetilde{M}(x,t)$ is clearly defined. Other methods for approximating nonlinear analytical solutions rely on pre-existing solutions, but we build the solution from scratch using only a minimal quantity of convergence-control parameters $C_i(i = 1, 2, ..., p)$ that are elements of the so-called optimal auxiliary functions.
The optimal auxiliary function method is a sequential approach that swiftly converges on the exact solution following the first iteration of the process. This method is based on the establishment and modification of auxiliary functions, as well as a simple mechanism for regulating the convergence of the solutions.

# 3 Numerical Examples and Applications

In this part of the article, we will look at the numerical solutions to several well-known nonlinear parabolic equations (the Benjamin-Bona-Mahony equation, Fisher's equation, the shock problem, and the Burger-Fisher equation) associated with initial conditions. In addition to the exact solution, graphical and numerical representations of all the numerical results and the absolute errors are provided here.
The absolute error refers to the discrepancy between the observed or estimated magnitude of a quantity and its real magnitude. The absolute error is determined by the following expression,

$$\text{Absolute Error (AE)} = |M(x,t) - \widetilde{M}(x,t)| \tag{16}$$

where $M(x,t)$ is the exact value and $\widetilde{M}(x,t)$ is the approximate value. The absolute error of the measurement shows how large the error actually is. In addition, the **Rate of Convergence** [43] can be defined as follows:

$$\mathcal{CR} = \frac{\log \frac{\epsilon_1}{\epsilon_2}}{\log \frac{t_1}{t_2}} \tag{17}$$

where $\epsilon_1$, $\epsilon_2$ are the absolute errors (AE) defined in Equation (16) for time step $t_1$ and $t_2$, respectively.

***Test Problem 1 :*** The Benjamin-Bona-Mahony equation models long waves in a nonlinear dispersive system. The equation is also called **Regularized Long-Wave Equation (RLWE)**. It is considered an improvement of the KdV equation. The equation covers the areas of surface waves of long wavelengths in liquids, acoustic-gravity waves in compressible fluids, hydromagnetic waves in a cold plasma, and acoustic waves in anharmonic crystals.
Let us consider the Benjamin-Bona-Mahony equation of the following type [28],

$$\left.\begin{aligned} &\frac{\partial M(x,t)}{\partial t} - \frac{\partial^3 M(x,t)}{\partial x^2 \partial t} + \frac{\partial M(x,t)}{\partial x} + M(x,t)\frac{\partial M(x,t)}{\partial x} = 0 \\ &M(x,0) = \operatorname{sech}^2\left(\frac{x}{4}\right) \end{aligned}\right\} \tag{18}$$

where $x \in [0, 0.07]$ and $t > 0$. The following provides an exact solution to the problem mentioned in (18)

$$M(x,t) = \operatorname{sech}^2\left(\frac{x}{4} - \frac{t}{3}\right)$$

The first-order approximate solution to the corresponding problem is given by,

$$\widetilde{M}(x,t) = \operatorname{sech}^2\left(\frac{x}{4}\right) + t\Bigg[C_1\left(\frac{1}{2}\operatorname{sech}^4\left(\frac{x}{4}\right)\tanh\left(\frac{x}{4}\right) + \frac{1}{2}\operatorname{sech}^6\left(\frac{x}{4}\right)\tanh\left(\frac{x}{4}\right)\right) + C_2\left(\frac{1}{2}\operatorname{sech}^6\left(\frac{x}{4}\right)\tanh\left(\frac{x}{4}\right) + \frac{1}{2}\operatorname{sech}^8\left(\frac{x}{4}\right)\tanh\left(\frac{x}{4}\right)\right) - C_3\operatorname{sech}^6\left(\frac{x}{4}\right) - C_4\operatorname{sech}^8\left(\frac{x}{4}\right)\Bigg] \tag{19}$$

The solution (19) can be written in the form of

$$\widetilde{M}(x,t) = M_0(x,t) + \sum_{j=1}^{n} tC_j\phi_j(x), \qquad n = 1,2,3,4 \tag{20}$$

where
$\phi_1(x) = \frac{1}{2}\operatorname{sech}^4\left(\frac{x}{4}\right)\tanh\left(\frac{x}{4}\right) + \frac{1}{2}\operatorname{sech}^6\left(\frac{x}{4}\right)\tanh\left(\frac{x}{4}\right),$
$\phi_2(x) = \frac{1}{2}\operatorname{sech}^6\left(\frac{x}{4}\right)\tanh\left(\frac{x}{4}\right) + \frac{1}{2}\operatorname{sech}^8\left(\frac{x}{4}\right)\tanh\left(\frac{x}{4}\right)$
$\phi_3(x) = -\operatorname{sech}^6\left(\frac{x}{4}\right),$
$\phi_4(x) = -\operatorname{sech}^8\left(\frac{x}{4}\right).$
The corresponding residual function can be represented as,

$$R(x,t) = \frac{\partial}{\partial t}\left(\widetilde{M}(x,t)\right) - \frac{\partial^3}{\partial x^2\partial t}\left(\widetilde{M}(x,t)\right) + \frac{\partial}{\partial x}\left(\widetilde{M}(x,t)\right) + \left(\widetilde{M}(x,t)\right)\frac{\partial}{\partial x}\left(\widetilde{M}(x,t)\right) \tag{21}$$

To obtain the values of auxiliary parameters' we set the residual equation as,

$$\int_0^{0.07} R(x,t)\phi_i(x)dx = 0$$

$$\text{or}\quad \sum_{j=1}^{n} C_j\int_0^{0.07}\left[\phi_j - \frac{\partial^2\phi_j}{\partial x^2} + t\frac{\partial\phi_j}{\partial x} + t\phi_j\left(\frac{\partial M_0}{\partial x} + \sum_{k=1}^{n} tC_k\frac{\partial\phi_k}{\partial x}\right) + M_0 t\frac{\partial\phi_j}{\partial x}\right]\phi_i(x)dx = \int_0^{0.07}\left[-\frac{\partial M_0}{\partial t} + \frac{\partial^3 M_0}{\partial x^2\partial t} - \frac{\partial M_0}{\partial x} - M_0\frac{\partial M_0}{\partial x}\right]\phi_i(x)dx$$

$$\text{or}\quad \sum_{j=1}^{n} C_j K_{ij} = F_i$$

where,

$$K_{ij} = \int_0^{0.07}\left[\phi_j - \frac{\partial^2\phi_j}{\partial x^2} + t\frac{\partial\phi_j}{\partial x} + t\phi_j\left(\frac{\partial M_0}{\partial x} + \sum_{k=1}^{n} tC_k\frac{\partial\phi_k}{\partial x}\right) + M_0 t\frac{\partial\phi_j}{\partial x}\right]\phi_i(x)dx$$

$$F_i = \int_0^{0.07}\left[-\frac{\partial M_0}{\partial t} + \frac{\partial^3 M_0}{\partial x^2\partial t} - \frac{\partial M_0}{\partial x} - M_0\frac{\partial M_0}{\partial x}\right]\phi_i(x)dx$$

The values of these coefficients are therefore inserted in the solution (19), yielding the final outcome. In the following table, the absolute errors obtained by different methods and our proposed method have been shown.

Table 1: Tabular Representations of the absolute error of (18) at different values of $x$

| $t$ | $x = 0.03$ | | | $x = 0.04$ | | |
|---|---|---|---|---|---|---|
| | **Absolute Error** | **Absolute Error [28]** | **Absolute Error[44]** | **Absolute Error** | **Absolute Error [28]** | **Absolute Error [44]** |
| 0.01 | $1.411273 \times 10^{-05}$ | $1.4104 \times 10^{-05}$ | $2.2664 \times 10^{-04}$ | $1.584927 \times 10^{-05}$ | $1.1584 \times 10^{-05}$ | $2.7703 \times 10^{-04}$ |
| 0.02 | $6.004954 \times 10^{-06}$ | $5.9887 \times 10^{-06}$ | $6.03525 \times 10^{-04}$ | $9.480540 \times 10^{-06}$ | $9.4689 \times 10^{-06}$ | $7.04304 \times 10^{-04}$ |
| 0.03 | $2.432482 \times 10^{-05}$ | $2.4349 \times 10^{-05}$ | $1.13601 \times 10^{-03}$ | $1.910947 \times 10^{-05}$ | $1.9126 \times 10^{-05}$ | $1.28165 \times 10^{-03}$ |
| 0.04 | $7.697610 \times 10^{-05}$ | $7.6908 \times 10^{-05}$ | $1.80786 \times 10^{-03}$ | $6.992144 \times 10^{-05}$ | $6.9944 \times 10^{-05}$ | $2.00908 \times 10^{-03}$ |
| 0.05 | $1.516464 \times 10^{-04}$ | $1.5168 \times 10^{-04}$ | $2.63254 \times 10^{-03}$ | $1.429545 \times 10^{-04}$ | $1.4298 \times 10^{-04}$ | $2.88653 \times 10^{-03}$ |

Table (1) has been used to represent the absolute errors at different time steps for values of $x = 0.03$ and $x = 0.04$. The collocation method was used to generate the values of the coefficients of solution (19) in reference [28]. In light of the information shown in Table (1), we are able to draw the conclusion that the solutions are pretty similar to one another. The test problem has been addressed using the optimal homotopy asymptotic methodology, as shown in reference [44]. The table shows that our findings are significantly better than those achieved by the preceding method.

***Test Problem 2 :*** Fisher's equation [45] is widely considered to be one of the most prominent examples of a nonlinear reaction-diffusion equation. The equation has been used in different aspects like flame propagation, chemical reactions, etc.
Let us now consider Fisher's equation with the initial condition,

$$\left.\begin{aligned} &\frac{\partial M(x,t)}{\partial t} = \frac{\partial^2 M(x,t)}{\partial x^2} + 6M(x,t)(1 - M(x,t)) \\ &M(x,0) = (1+e^x)^{-2} \end{aligned}\right\} \tag{22}$$

where $x \in [0,1]$ and $t > 0$. The following equation provides an exact solution to the aforementioned problem (22),

$$M(x,t) = (1+e^{x-5t})^{-2}.$$

The first-order approximate solution to the corresponding problem is given by,

$$\begin{aligned} \widetilde{M}(x,t) = (1+e^x)^{-2} - t\Big[C_1(1+e^x)^{-2} + C_2\Big((1+e^x)^{-2}\Big)^2 \\ + C_3\Big((1+e^x)^{-2}\Big)^3 + C_4\Big((1+e^x)^{-2}\Big)^4\Big] \end{aligned} \tag{23}$$

where
$\phi_1(x) = -(1+e^x)^{-2}$,
$\phi_2(x) = -\Big((1+e^x)^{-2}\Big)^2$,
$\phi_3(x) = -\Big((1+e^x)^{-2}\Big)^3$,
$\phi_4(x) = -\Big((1+e^x)^{-2}\Big)^4$.
The residual function $R(x,t)$ can be defined as

$$R(x,t) = \frac{\partial \widetilde{M}(x,t)}{\partial t} - \frac{\partial^2 \widetilde{M}(x,t)}{\partial x^2} - 6\widetilde{M}(x,t)(1 - \widetilde{M}(x,t)) \tag{24}$$

and the residual equation can be written as

$$\int_0^1 R(x,t)\phi_i(x)dx = 0$$

$$\text{or} \sum_{j=1}^{n} C_j \int_0^1 \left[\phi_j - t\frac{\partial^2 \phi_j}{\partial x^2} + 6tM_0\phi_j - 6t\phi_j + 6tM_0\phi_j + \left(\sum_{k=1}^{n} 6t^2 C_k\phi_k\right)\phi_j\right]\phi_i(x)dx$$
$$= \int_0^1 \left[-\frac{\partial M_0}{\partial t} + \frac{\partial^2 M_0}{\partial x^2} + 6M_0 - 6M_0^2\right]\phi_i(x)dx$$

$$\text{or} \quad \sum_{j=1}^{n} C_j K_{ij} = F_i$$

where

$$K_{ij} = \int_0^1 \left[\phi_j - t\frac{\partial^2 \phi_j}{\partial x^2} + 6tM_0\phi_j - 6t\phi_j + 6tM_0\phi_j + \left(\sum_{k=1}^{n} 6t^2 C_k\phi_k\right)\phi_j\right]\phi_i(x)dx$$
$$F_i = \int_0^1 \left[-\frac{\partial M_0}{\partial t} + \frac{\partial^2 M_0}{\partial x^2} + 6M_0 - 6M_0^2\right]\phi_i(x)dx$$

We have evaluated the values of these coefficients by solving the system of equations. The values are therefore inserted in the solution (23), yielding the final outcome.

**Table 2:** Tabular Representations of approximate solutions and absolute errors of (22) at different time steps

| $x$ | $t = 0.001$ | | | $t = 0.01$ | | |
|---|---|---|---|---|---|---|
| | **Approximate Solution** | **Absolute Error** | **Absolute Error[45]** | **Approximate Solution** | **Absolute Error** | **Absolute Error [45]** |
| 0.0 | 0.25125226 | $7.00945 \times 10^{-07}$ | $5.0 \times 10^{-04}$ | 0.26274282 | $8.92336 \times 10^{-05}$ | $4.7 \times 10^{-03}$ |
| 0.1 | 0.22683291 | $1.84484 \times 10^{-06}$ | $1.0 \times 10^{-04}$ | 0.23780423 | $1.45442 \times 10^{-04}$ | $7.0 \times 10^{-04}$ |
| 0.2 | 0.20376738 | $1.90543 \times 10^{-06}$ | $2.0 \times 10^{-04}$ | 0.21414333 | $1.72178 \times 10^{-04}$ | $2.3 \times 10^{-03}$ |
| 0.3 | 0.18214311 | $1.75098 \times 10^{-06}$ | $2.0 \times 10^{-04}$ | 0.19187472 | $1.85299 \times 10^{-04}$ | $2.2 \times 10^{-03}$ |
| 0.4 | 0.16201943 | $1.72015 \times 10^{-06}$ | $2.0 \times 10^{-04}$ | 0.17107753 | $1.92511 \times 10^{-04}$ | $2.4 \times 10^{-03}$ |
| 0.5 | 0.14342796 | $1.84292 \times 10^{-06}$ | $2.0 \times 10^{-04}$ | 0.15179831 | $1.96502 \times 10^{-04}$ | $2.5 \times 10^{-03}$ |
| 0.6 | 0.12637404 | $2.00938 \times 10^{-06}$ | $3.0 \times 10^{-04}$ | 0.13405416 | $1.97392 \times 10^{-04}$ | $2.5 \times 10^{-03}$ |
| 0.7 | 0.11083896 | $2.08225 \times 10^{-06}$ | $2.0 \times 10^{-04}$ | 0.11783625 | $1.94426 \times 10^{-04}$ | $2.5 \times 10^{-03}$ |
| 0.8 | 0.09678273 | $1.95933 \times 10^{-06}$ | $2.0 \times 10^{-04}$ | 0.10311327 | $1.86959 \times 10^{-04}$ | $2.2 \times 10^{-03}$ |
| 0.9 | 0.08414747 | $1.59713 \times 10^{-06}$ | $2.0 \times 10^{-04}$ | 0.08983494 | $1.74906 \times 10^{-04}$ | $2.6 \times 10^{-03}$ |
| 1.0 | 0.07286085 | $1.00891 \times 10^{-06}$ | $1.0 \times 10^{-04}$ | 0.07793554 | $1.58800 \times 10^{-04}$ | $6.0 \times 10^{-03}$ |

Table (2) has been yielded to represent the approximate solutions and absolute errors of (22) at different time steps. The table has also included the reference absolute errors from the previously published literature. Based on Table (2), it can be shown that our suggested methodology provides a high level of precision across the domain and at various time levels. Comparing the features of the approximate solution with those of the exact solution reveals a high degree of similarity.
The following table shows the convergence rate of the test problem using our proposed method for different time steps.

**Table 3:** Convergence Rate ($\mathcal{CR}$) using the present approach for Test Problem 2

| $t$ | **MAE** | **Convergence Rate** |
|---|---|---|
| 0.001 | $2.08 \times 10^{-06}$ | |
| 0.002 | $7.890 \times 10^{-06}$ | 1.9234 |
| 0.003 | $1.758 \times 10^{-05}$ | 1.9759 |
| 0.004 | $3.122 \times 10^{-05}$ | 1.9963 |
| 0.005 | $4.882 \times 10^{-05}$ | 2.0035 |

Figure (1) illustrates the required results derived from equation (22) employing the developed algorithm of our provided method for different time steps. Given how similar the approximate and exact solutions are, distinguishing between them through these diagrams is a challenging task.

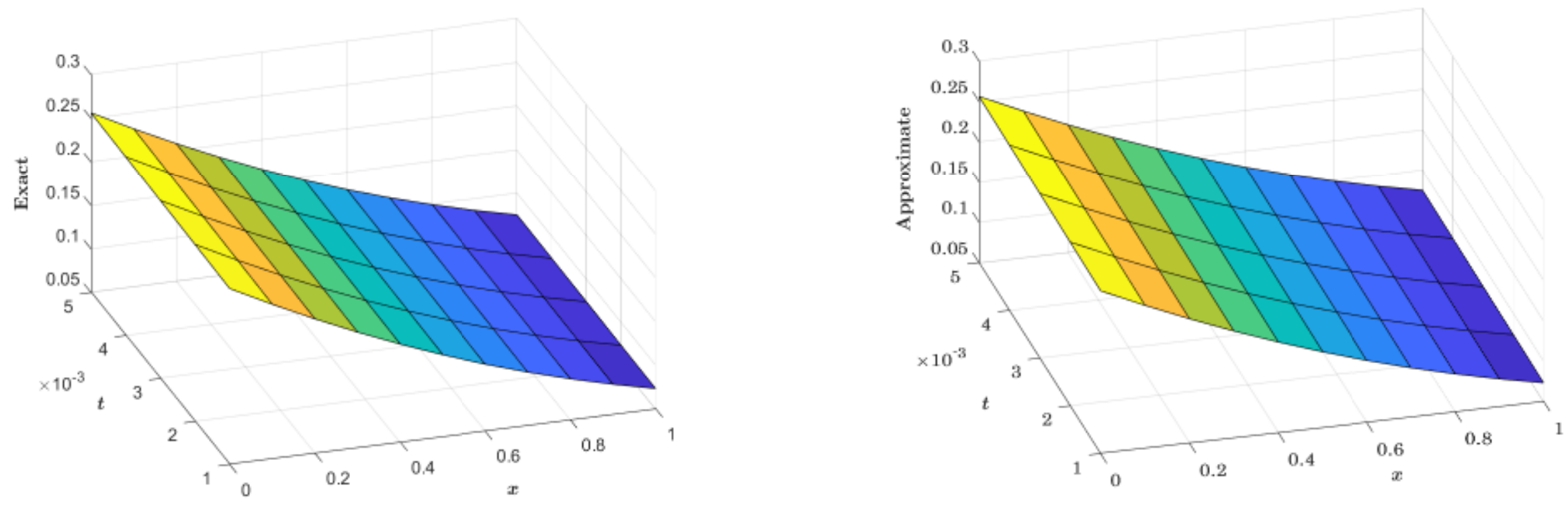

**(a)** **(b)**

**Figure 1:** Exact (on left) and Approximate (on right) solution of (22) at different time steps

Figure (2) displays the error graphs that show the absolute discrepancy between the numerical and exact solutions of (22). An acceptable level of inaccuracy is provided by the absolute error map for the employment of the OAFM in illustrating numerical results.

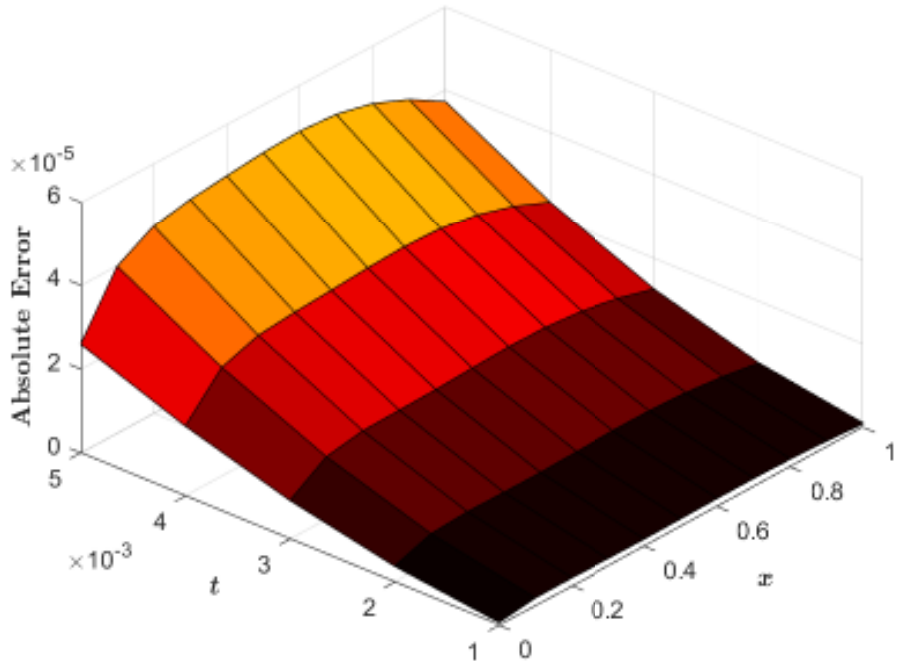

**Figure 2:** The schematic 3D representation of absolute error of (22) for different time steps

***Test Problem 3:*** The shock wave is a prominent type of equation used in science and technology. The shock wave has interesting applications in a variety of areas, such as medicine, biological sciences, material processing, manufacturing, and microelectronic industries. Let us now take into consideration the uniformly propagating shock problem [46],

$$\left.\begin{aligned} &\frac{\partial M(x,t)}{\partial t} = \frac{1}{Re}\frac{\partial^2 M(x,t)}{\partial x^2} - M(x,t)\frac{\partial M(x,t)}{\partial x} \\ &M(x,0) = \frac{x-4}{x-2} \end{aligned}\right\} \tag{25}$$

where $x \in [-1,1]$, $t > 0$ and $Re$ is known as the Reynolds number. Here $Re = 1$.
The following provides an exact solution to the aforementioned problem (25),

$$M(x,t) = 1 - \frac{2}{x-t-2}$$

The first-order approximate solution of (25) is obtained as,

$$\widetilde{M}(x,t) = \frac{x-4}{x-2} - t\left[C_1\frac{x-4}{x-2} + C_2\left(\frac{x-4}{x-2}\right)^2 + C_3\left(\frac{x-4}{x-2}\right)^3 + C_4\left(\frac{x-4}{x-2}\right)^4\right] \tag{26}$$

where

$\phi_1(x) = -\left(\frac{x-4}{x-2}\right),$

$\phi_2(x) = -\left(\frac{x-4}{x-2}\right)^2,$

$\phi_3(x) = -\left(\frac{x-4}{x-2}\right)^3,$

$\phi_4(x) = -\left(\frac{x-4}{x-2}\right)^4.$

The residual function $R(x,t)$ can be defined as

$$R(x,t) = \frac{\partial \widetilde{M}(x,t)}{\partial t} - \frac{1}{Re}\frac{\partial^2 \widetilde{M}(x,t)}{\partial x^2} + \widetilde{M}(x,t)\frac{\partial \widetilde{M}(x,t)}{\partial x} \tag{27}$$

and the residual equation can be written as

$$\int_{-1}^{1} R(x,t)\phi_i(x)dx = 0$$

$$\text{or} \sum_{j=1}^{n} C_j \int_{-1}^{1}\left[\phi_j - \frac{t}{Re}\frac{\partial^2\phi_j}{\partial x^2} + tM_0\frac{\partial \phi_j}{\partial x} + t\phi_j\frac{\partial M_0}{\partial x} + \left(\sum_{k=1}^{n} t^2 C_k\phi_k\right)\phi_j\right]\phi_i(x)dx = \int_{-1}^{1}\left[-\frac{\partial M_0}{\partial t} + \frac{1}{Re}\frac{\partial^2 M_0}{\partial x^2} + M_0\frac{\partial M_0}{\partial x}\right]\phi_i(x)dx$$

$$\text{or} \quad \sum_{j=1}^{n} C_j K_{ij} = F_i$$

where

$$K_{ij} = \int_{-1}^{1}\left[\phi_j - \frac{t}{Re}\frac{\partial^2\phi_j}{\partial x^2} + tM_0\frac{\partial \phi_j}{\partial x} + t\phi_j\frac{\partial M_0}{\partial x} + \left(\sum_{k=1}^{n} t^2 C_k\phi_k\right)\phi_j\right]\phi_i(x)dx$$

$$F_i = \int_{-1}^{1}\left[-\frac{\partial M_0}{\partial t} + \frac{1}{Re}\frac{\partial^2 M_0}{\partial x^2} + M_0\frac{\partial M_0}{\partial x}\right]\phi_i(x)dx$$

In order to guarantee the efficacy and dependability of this process, a comparison will be made between the approximate numerical solution and the exact solution in Table (4).

**Table 4:** Tabular Representations of approximate solutions and absolute errors of (25) at different time steps

| $x$ | $t = 0.01$ | | $t = 0.02$ | | $t = 0.03$ | |
|---|---|---|---|---|---|---|
| | **Approximate** | **Absolute Error** | **Approximate** | **Absolute Error** | **Approximate** | **Absolute Error** |
| -1.0 | 1.66446466 | $1.28377 \times 10^{-05}$ | 1.66228964 | $7.00945 \times 10^{-05}$ | 1.66013535 | $6.93521 \times 10^{-05}$ |
| -0.8 | 1.71175467 | $1.09049 \times 10^{-05}$ | 1.70925886 | $1.84484 \times 10^{-05}$ | 1.70679439 | $8.06126 \times 10^{-05}$ |
| -0.6 | 1.76629410 | $1.05769 \times 10^{-05}$ | 1.76340273 | $1.90543 \times 10^{-05}$ | 1.76055458 | $9.83131 \times 10^{-04}$ |
| -0.4 | 1.82988779 | $1.22779 \times 10^{-05}$ | 1.82650026 | $1.75098 \times 10^{-05}$ | 1.82316979 | $1.24524 \times 10^{-04}$ |
| -0.2 | 1.90499387 | $1.64933 \times 10^{-05}$ | 1.90097158 | $1.72015 \times 10^{-05}$ | 1.89702342 | $1.62435 \times 10^{-04}$ |
| 0.0 | 1.99504863 | $2.37565 \times 10^{-05}$ | 1.99019534 | $1.84292 \times 10^{-05}$ | 1.98543925 | $2.17583 \times 10^{-04}$ |
| 0.2 | 2.10500705 | $3.46761 \times 10^{-05}$ | 2.09903586 | $2.00938 \times 10^{-04}$ | 2.09319728 | $3.01108 \times 10^{-04}$ |
| 0.4 | 2.24228626 | $5.02384 \times 10^{-05}$ | 2.23476196 | $2.08225 \times 10^{-04}$ | 2.22743305 | $4.39187 \times 10^{-04}$ |
| 0.6 | 2.41851329 | $7.35720 \times 10^{-05}$ | 2.40874659 | $1.95933 \times 10^{-04}$ | 2.39930454 | $7.03143 \times 10^{-04}$ |
| 0.8 | 2.65301155 | $1.18993 \times 10^{-04}$ | 2.63985645 | $1.59713 \times 10^{-04}$ | 2.62733328 | $1.31702 \times 10^{-03}$ |
| 1.0 | 2.98045785 | $2.59828 \times 10^{-04}$ | 2.96191353 | $1.00891 \times 10^{-03}$ | 2.94484991 | $3.10234 \times 10^{-03}$ |

According to Table (4), it is evident that our aforementioned method guarantees high precision across a variety of time scales.

Figure (3) provides the 3D visual representations of the exact and approximate results of (25) for different time steps. As can be seen in Figure (3), it is difficult to differentiate between exact and approximate solutions, which paves the way for the acceptance of numerical representations derived through the implementation of the optimal auxiliary function method.

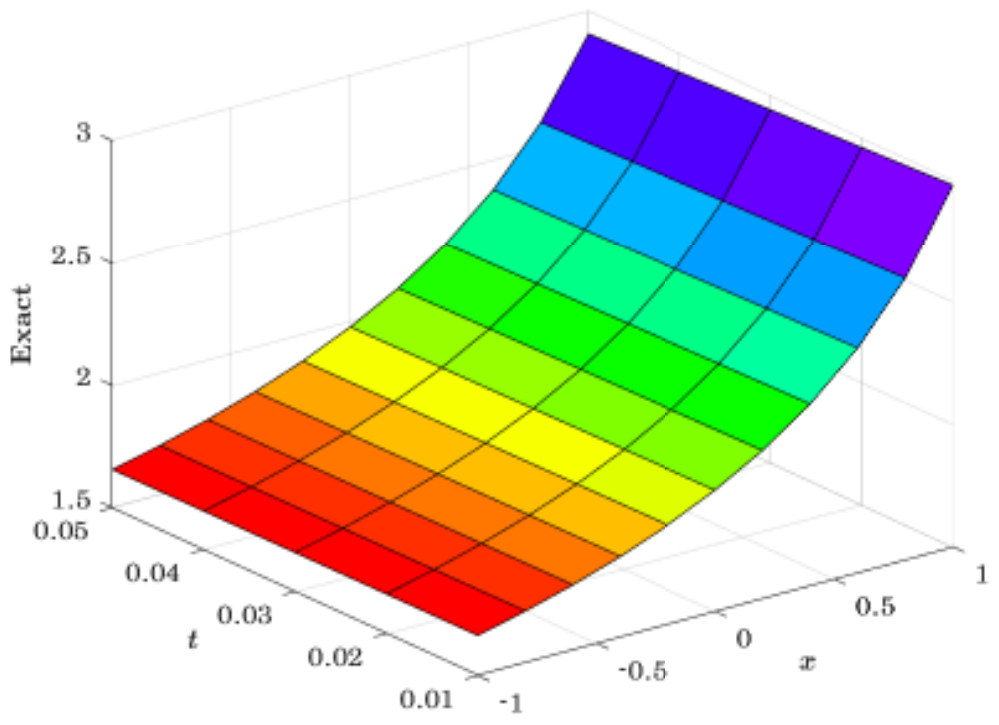

**(a)**

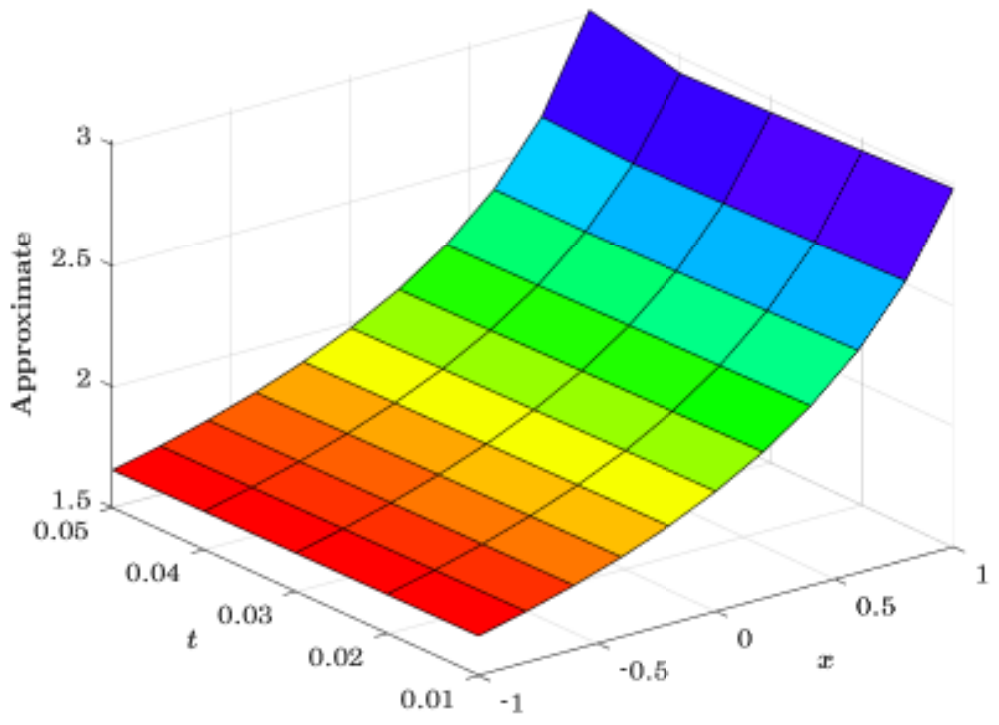

**(b)**

**Figure 3:** Exact (on left) and Approximate (on right) solution with of (25) at different time steps

Figure (4) exhibits a visual representation of the absolute errors for different time steps and for different values of $x$. The figure shows that the absolute errors are quite negligible. In addition, the figure verifies the solution's reliability and acceptance of the recommended approach.

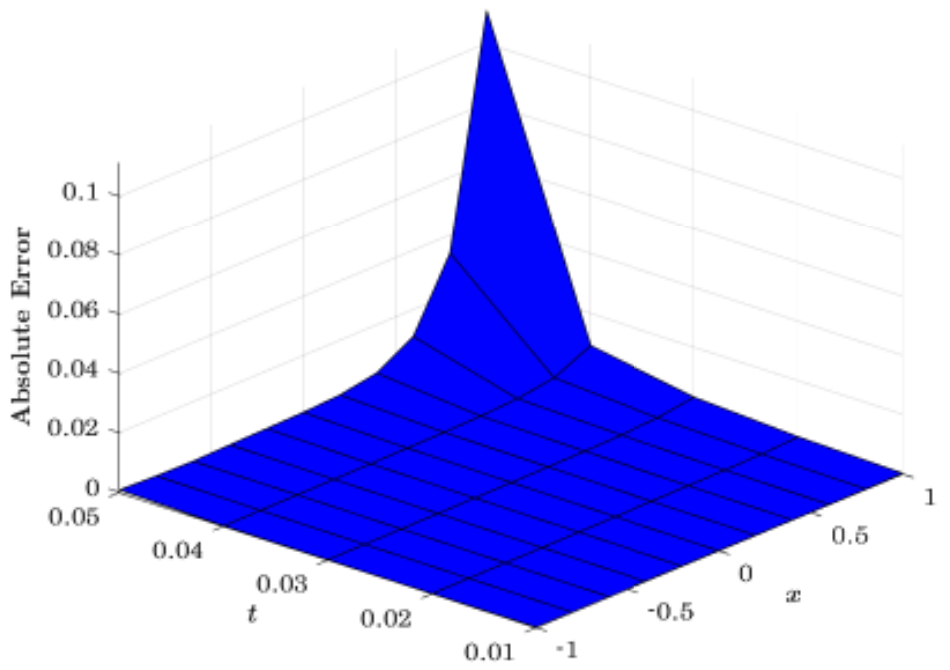

**Figure 4:** The schematic representation of absolute error of (25) for various time steps

***Test Problem 4:*** The Burger-Fisher equation illustrates a conventional concept for modeling the interplay between the reaction process, the convection impact, and diffusive mobility. The equation is important to the study of a wide variety of subfields within mathematics and physics, including finance, gas dynamics, and traffic patterns. Let us consider the Burger-Fisher equation [47] with the initial condition,

$$\left.\begin{aligned} &\frac{\partial M(x,t)}{\partial t} + \alpha M(x,t)^{\omega}\frac{\partial M(x,t)}{\partial x} - \frac{\partial^2 M(x,t)}{\partial x^2} = \beta M(x,t)(1 - M(x,t)^{\omega}) \\ &M(x,0) = \left\{\frac{1}{2} + \frac{1}{2}\tanh\left(\frac{-\alpha\omega}{2(\omega+1)}x\right)\right\}^{\frac{1}{\omega}} \end{aligned}\right\} \tag{28}$$

where $x \in [0,1]$ and $t > 0$. We choose $\alpha = 1$, $\beta = 1$, and $\omega = 1$ to get our approximate solution. The following provides an exact solution to the aforementioned problem (28),

$$M(x,t) = \left(\frac{1}{2} + \frac{1}{2}\tanh\left[\frac{-\alpha\omega}{2(\omega+1)}\left(x - \left(\frac{\alpha}{\omega+1} + \frac{\beta(\omega+1)}{\alpha}\right)t\right)\right]\right)^{\frac{1}{\omega}}.$$

The first-order approximate solution of (28) is obtained as,

$$\begin{aligned} \widetilde{M}(x,t) = \left(\frac{1}{2} + \frac{1}{2}\tanh\left(\frac{-x}{4}\right)\right) - t\Bigg[C_1\left(\frac{1}{2} + \frac{1}{2}\tanh\left(\frac{-x}{4}\right)\right) + C_2\left(\frac{1}{2} + \frac{1}{2}\tanh\left(\frac{-x}{4}\right)\right)^2 \\ + C_3\left(\frac{1}{2} + \frac{1}{2}\tanh\left(\frac{-x}{4}\right)\right)^3 + C_4\left(\frac{1}{2} + \frac{1}{2}\tanh\left(\frac{-x}{4}\right)\right)^4\Bigg] \end{aligned} \tag{29}$$

The solution (29) can be written in the form of

$$\widetilde{M}(x,t) = M_0(x,t) + \sum_{j=1}^{n} tC_j\phi_j(x), \qquad n = 1,2,3,4 \tag{30}$$

where

$\phi_1(x) = -\left(\frac{1}{2} + \frac{1}{2}\tanh\left(\frac{-x}{4}\right)\right)$,

$\phi_2(x) = -\left(\frac{1}{2} + \frac{1}{2}\tanh\left(\frac{-x}{4}\right)\right)^2$,

$\phi_3(x) = -\left(\frac{1}{2} + \frac{1}{2}\tanh\left(\frac{-x}{4}\right)\right)^3$,

$\phi_4(x) = -\left(\frac{1}{2} + \frac{1}{2}\tanh\left(\frac{-x}{4}\right)\right)^4$.

The corresponding residual function can be represented as,

$$R(x,t) = \frac{\partial \widetilde{M}}{\partial t} + \widetilde{M}\frac{\partial \widetilde{M}}{\partial x} - \frac{\partial^2 \widetilde{M}}{\partial x^2} - \widetilde{M}(1-\widetilde{M}) \tag{31}$$

To obtain the values of auxiliary parameters' we set the residual equation as,

$$\begin{aligned}
&\int_0^1 R(x,t)\phi_i(x)dx = 0 \\
\text{or} \quad &\sum_{j=1}^{n} C_j \int_0^1 \Bigg[\phi_j - t\phi_j\frac{\partial M_0}{\partial x} + M_0 t\frac{\partial \phi_j}{\partial x} + \Big(\sum_{k=1}^{n} tC_k\frac{\partial \phi_k}{\partial x}\Big)\phi_j - t\frac{\partial^2 \phi_j}{\partial x^2} - t\phi_j(1-M_0) + M_0 t\phi_j \\
&\quad + t^2\phi_j\Big(\sum_{k=1}^{n} C_k\phi_k\Big)\Bigg]\phi_i(x)dx = \int_0^1 \Bigg[-\frac{\partial M_0}{\partial t} - M_0\frac{\partial M_0}{\partial x} + \frac{\partial^2 M_0}{\partial x^2} + M_0(1-M_0)\Bigg]\phi_i(x)dx \\
\text{or} \quad &\sum_{j=1}^{n} C_j K_{ij} = F_i
\end{aligned}$$

where,

$$\begin{aligned}
K_{ij} &= \int_0^1 \Bigg[\phi_j - t\phi_j\frac{\partial M_0}{\partial x} + M_0 t\frac{\partial \phi_j}{\partial x} + \Big(\sum_{k=1}^{n} tC_k\frac{\partial \phi_k}{\partial x}\Big)\phi_j - t\frac{\partial^2 \phi_j}{\partial x^2} - t\phi_j(1-M_0) + M_0 t\phi_j \\
&\qquad + t^2\phi_j(\sum_{k=1}^{n} C_k\phi_k)\Bigg]\phi_i(x)dx \\
F_i &= \int_0^1 \Bigg[-\frac{\partial M_0}{\partial t} - M_0\frac{\partial M_0}{\partial x} + \frac{\partial^2 M_0}{\partial x^2} + M_0(1-M_0)\Bigg]\phi_i(x)dx
\end{aligned}$$

Then the values of the coefficients $C_1$, $C_2$, $C_3$, and $C_4$ are then obtained by solving the system of equations. Then we substitute these coefficients' values in solution (29) and the final approximated result of (28) is obtained.

**Table 5:** Tabular Representations of approximate solutions and absolute errors of (28) at different time steps

| $x$ | $t = 0.01$ | | $t = 0.05$ | | $t = 0.10$ | |
|---|---|---|---|---|---|---|
| | **Approximate** | **Absolute Error** | **Approximate** | **Absolute Error** | **Approximate** | **Absolute Error** |
| 0.0 | 0.50311707 | $7.88760 \times 10^{-06}$ | 0.51541589 | $2.04019 \times 10^{-04}$ | 0.53036447 | $8.44902 \times 10^{-04}$ |
| 0.1 | 0.49061870 | $7.39628 \times 10^{-06}$ | 0.50293257 | $1.92383 \times 10^{-04}$ | 0.51793963 | $8.01589 \times 10^{-04}$ |
| 0.2 | 0.47813206 | $6.88732 \times 10^{-06}$ | 0.49044584 | $1.80254 \times 10^{-04}$ | 0.50549354 | $7.56136 \times 10^{-04}$ |
| 0.3 | 0.46567269 | $6.36334 \times 10^{-06}$ | 0.47797125 | $1.67690 \times 10^{-04}$ | 0.49304159 | $7.08728 \times 10^{-04}$ |
| 0.4 | 0.45325602 | $5.82707 \times 10^{-06}$ | 0.46552430 | $1.54751 \times 10^{-04}$ | 0.48059921 | $6.59576 \times 10^{-04}$ |
| 0.5 | 0.44089725 | $5.28131 \times 10^{-06}$ | 0.45312034 | $1.41505 \times 10^{-04}$ | 0.46818172 | $6.08911 \times 10^{-04}$ |
| 0.6 | 0.42861129 | $4.72898 \times 10^{-06}$ | 0.44077450 | $1.28020 \times 10^{-04}$ | 0.45580433 | $5.56982 \times 10^{-04}$ |
| 0.7 | 0.41641271 | $4.17297 \times 10^{-06}$ | 0.42850165 | $1.14367 \times 10^{-04}$ | 0.44348206 | $5.04051 \times 10^{-04}$ |
| 0.8 | 0.40431565 | $3.61619 \times 10^{-06}$ | 0.41631626 | $1.00616 \times 10^{-04}$ | 0.43122962 | $4.50392 \times 10^{-04}$ |
| 0.9 | 0.39233377 | $3.06152 \times 10^{-06}$ | 0.40423242 | $8.68417 \times 10^{-05}$ | 0.41906141 | $3.96286 \times 10^{-04}$ |
| 1.0 | 0.38048017 | $2.51176 \times 10^{-06}$ | 0.39226371 | $7.31134 \times 10^{-05}$ | 0.40699138 | $3.42017 \times 10^{-04}$ |

Table (5) provides approximate results and absolute errors of the corresponding problem at different time levels. The table assures that error terms are quite minimal, and enhances the chances of the suggested scheme being approved.

The following table shows the convergence rate of the test problem using our proposed method for different time steps.

**Table 6:** Convergence Rate ($\mathcal{CR}$) using the present approach for Test Problem 4

| $t$ | MAE | Convergence Rate |
|---|---|---|
| 0.01 | $7.8876 \times 10^{-06}$ | |
| 0.02 | $3.1840 \times 10^{-05}$ | 2.0132 |
| 0.03 | $7.2265 \times 10^{-05}$ | 2.0214 |
| 0.04 | $1.2954 \times 10^{-04}$ | 2.0288 |
| 0.05 | $2.0402 \times 10^{-04}$ | 2.0356 |

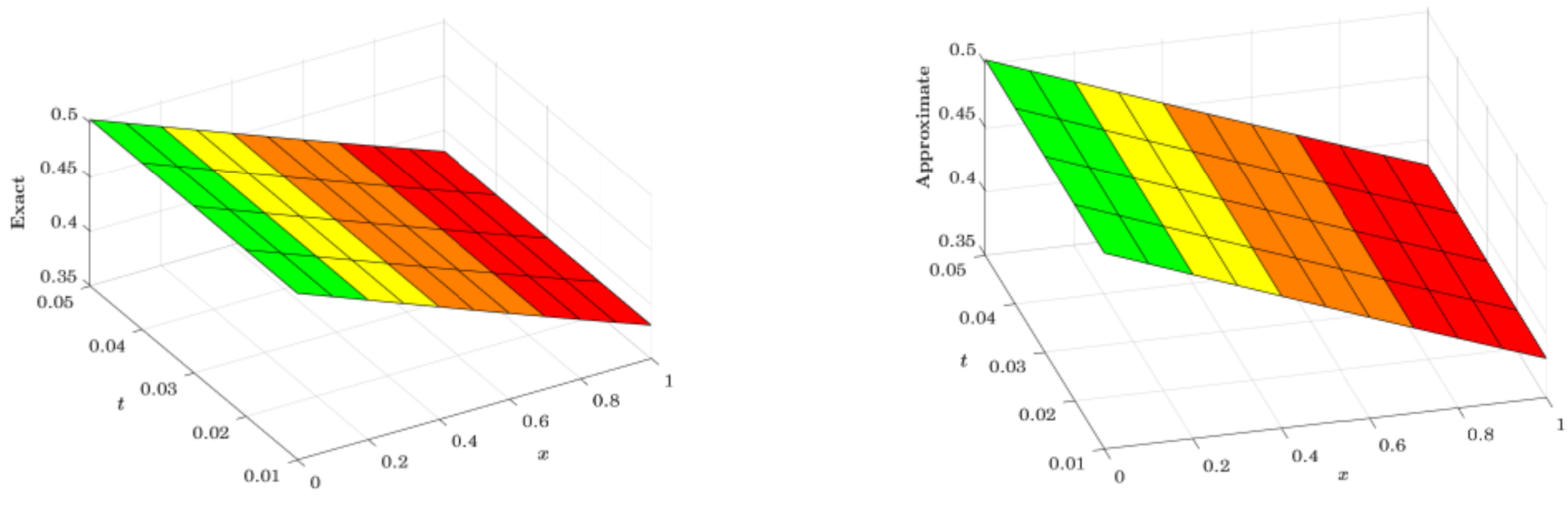

(a) (b)

**Figure 5:** Exact (on left) and Approximate (on right) solution with of (28) at different time steps

Figure (5) is deployed to provide pictorial representations of the exact and approximate data obtained by the present methodology at different time steps. The figure represents a good agreement between the exact and approximate data.

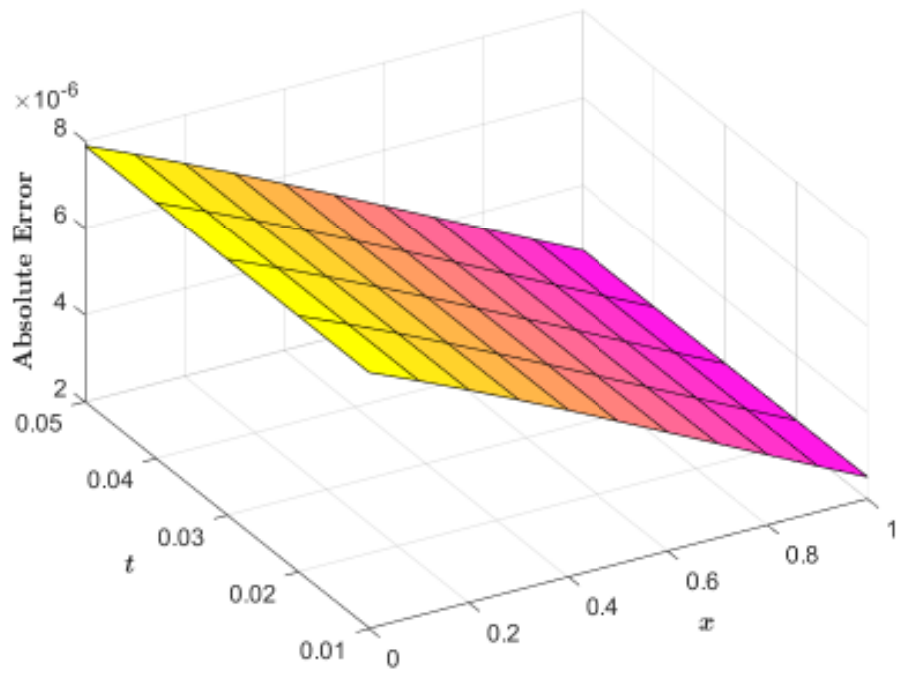

**Figure 6:** The schematic representation of absolute error of (28) for different time steps

Figure (6) exhibits a visual representation of the absolute errors for various nodes at different time steps. The figure ensures that the algorithm is reliable enough.

# Conclusion

In this research, we have obtained numerical approximations of nonlinear parabolic PDEs with the help of the initial conditions through the application of the optimal auxiliary function method. This procedure has been put in place to determine the initial and first estimations of the parabolic PDEs. Then, by adding these parts, we have derived the first-order approximate solutions. We have also employed the well-known Galerkin approach to calculate the coefficients of initial estimation. Then, the technique has been applied to numerically solve a number of well-known parabolic PDEs, with the approximation results being compared to the exact ones. All of the approximate results, together with their graphical and tabular representations, as well as their absolute errors, are provided. In light of these facts, it is clear that the approach that was suggested is efficient, reliable, and acceptable due to the fact that it is easy to implement. In the future, the proposed method can be applied to 2D and 3D nonlinear parabolic PDEs.

## Acknowledgement

The first author is grateful to the National Science & Technology (NST), Ministry of Science & Technology, Govt. of the People's Republic of Bangladesh, for granting the "NST Fellowship" partially during the period of research work.